
\emergencystretch=5mm
\def\pp{ {\bf P} }
\def\rr{ {\bf R} }
\def\cc{ {\bf C} }
\def\Ee{ {\cal E} }
\def\Oo{ {\cal O} }
\def\Ll{ {\cal L} }
\def\ext{\mathop{\rm Ext}\nolimits^1}
\def\inPn{\subseteq\pp^n}
\def\map#1#2#3{#1\colon#2\rightarrow #3}
\def\stackrel#1#2{\mathrel{\mathop{#2}\limits^{#1}}}
\def\triple#1#2#3#4#5{0\rightarrow#1
     \stackrel{#2}{\rightarrow}#3%
     \stackrel{#4}{\rightarrow}#5%
     \rightarrow 0}
\def\frac#1#2{{#1\over#2}}
\def\vertmap#1{\Big\downarrow\rlap
                            {$\vcenter{\hbox{$
                \displaystyle#1$}}$}
}
\def \proj {\mathop{\rm Proj}\nolimits}
\def \sym {\mathop{\rm Sym}\nolimits}
\def \dim {\mathop{\rm dim}\nolimits}
\def \Aut {\mathop{\rm Aut}\nolimits}
\def \Ker {\mathop{\rm Ker}\nolimits}
\def \diag {\mathop{\rm diag}\nolimits}
\def\Im{ {\rm Im} }
%
%
\leftline{AMS subject classification: 14D05}\bigskip
\centerline{S.L'vovsky\footnote{*}{Partially supported
by the AMS grant for mathematicians in former Soviet Union.}}\bigskip

\centerline{ON CURVES AND SURFACES WITH}
\centerline{PROJECTIVELY EQUIVALENT}
\centerline{HYPERPLANE SECTIONS}
\bigskip

\centerline{\bf Introduction}\medskip
In this paper we are concerned with the following question: how
to describe the projective varieties such that almost all their
hyperplane sections are projectively equivalent? We give the
complete answer for curves and a partial one for smooth
surfaces (in characteristic~0 both).

The question we are interested in was considered, for the case
of surfaces, in 1925 by Guido Fubini and Gino Fano ([6,4,5]).
The final results are contained in~[5]. Fano gives the complete
list of surfaces with projectively equivalent hyperplane
sections (and arbitrary singularities); we consider only smooth
surfaces, and our list is apparently superfluous: according to
Fano, some of the surfaces therein should not have projectively
equivalent hyperplane sections, but I did not manage to prove
it, nor to follow the argument in~[5]. For the case of curves
in characteristic 0, our result is complete.

Nowadays this problem was considered by Edoardo Ballico~[1] in
arbitrary characteristic. Our method differs from that of~[1].
For the case of curves our result is in accord with~[1], for
the case of surfaces in characteristic 0 our result strengthens
Proposition~5.2 of~[1].

In the appendix we prove a result concerning connections
between projective equivalence of hyperplane sections,
finiteness of monodromy group and the adjunction properties of
a variety.

When the first draft of this paper was finished, I learned that
Rita Pardini~[9] had proved Fano results from~[5] in full.

\bigskip
\centerline{\bf Notation and conventions}\medskip
Throughout the paper, the base field will be the field $\cc$ of
complex numbers. If  $\Ee$ is a  locally free sheaf  over  $X$,
then $\pp(\Ee)=\proj\sym(\Ee)$.

If $p\colon X\rightarrow  Y$ is a nonramified covering, we will
say that it is {\it split\/} if each connected component of $X$
is mapped isomorphically on $Y$.

We will say that a surface $X\subseteq\pp^n$ is a {\it
scroll\/} if
$$(X,\Oo_X(1))\cong
(\pp(\Ee),\Oo_{\pp(\Ee)|C}(1))
$$
for a smooth curve $C$ and a
rank 2 locally free sheaf $\Ee$.

If $(\pp^n)^\ast$ is the dual to projective space $\pp^n$ and
$\alpha\in (\pp^n)^\ast$, we denote by $H_\alpha\subset\pp^n$
the hyperplane corresponding to $\alpha $.

We say that a projective variety $X\subseteq\pp^n$ is {\it
linearly normal\/} if the linear system of its hyperplane
sections is complete.

If $X\subseteq\pp^n$ is a smooth projective variety of
dimension $d$, the monodromy group of its hyperplane section is
the monodromy group acting on $H^{d-1}(Y,\rr)$ as its smooth
hyperplane section $Y$ varies (cf.~[3]).

\bigskip
\centerline{\bf Statement of results}\medskip
Let $X\subseteq\pp^n$ be a projective variety. We say that $X$
satisfies {\it the FF condition\/} (named so after G.Fubini and
G.Fano) if almost all hyperplane sections of $X$ are
projectively equivalent.

\proclaim Proposition 0{.}1.
If $X$ is an irreducible curve not contained in  a  hyperplane,
then FF condition is satisfied if and only if $\deg X\le n+1$.

\proclaim Proposition 0{.}2.
If $X\subseteq\pp^n$ is a smooth irreducible surface satisfying
the FF condition, then $X$ is either a rational scroll, or a
Veronese surface $v_2(\pp^2)\subset\pp^5$, or its isomorphic
projection, or a non-linearly-normal scroll with elliptic base.

According to Fano~[5], of all the surfaces listed in the above
proposition, only linearly normal rational scrolls and
$v_2(\pp^2)\subset\pp^5$ satisfy the FF condition.

     Here is an amusing corollary to Proposition~0.2.

\proclaim Proposition 0{.}3.
If the surface $X\subseteq\pp^n$ is a scroll with base of genus
$>2$ or a linearly normal scroll with elliptic base, then
almost all hyperplane sections of $X$ are not linearly normal.

\bigskip
\centerline{\bf Acknowledgements}\medskip
I am grateful to F.L.Zak for attracting my attention to the
papers~[6,4,5] and for numerous helpful discussions. I would
like to thank E.Ballico and Rita Pardini for helpful
correspondence, and S.Tabachnikov for helpful discussions.

\bigskip
\centerline{\bf 1. Preliminaries; the FF condition and monodromy}
\medskip
Let $X\inPn$ be a projective variety. For
$\alpha,\beta\in(\pp^n)^\ast$ denote by $\Phi_{\alpha\beta}$
the set of projective isomorphisms
$\map{\varphi}{H_\alpha}{H_\beta}$ such that $\varphi(X\cap
H_\alpha)=X\cap H_\beta$.

\proclaim Proposition 1{.}1.
Assume that $X$ satisfies the FF condition and $X$ is smooth or
$\dim X=1$; if we denote by $Y=X\cap\pp^n$ the generic
hyperplane section, then the action of the monodromy group on
$H^\bullet(Y)$ is induced by the action of a subgroup
$G\subseteq \{g\in \Aut(\pp^n): gY=Y\}$.

{\bf Proof.} The condition FF implies that there exists a
Zariski open subset $U\subseteq (\pp^n)^\ast$ such that, for
$\alpha,\beta\in U$ we have $\Phi_{\alpha\beta}\ne\emptyset$
and $H_\alpha$ is transversal to $X$. Consider a fiber space
$\Phi$ over $U\times U$ such that $\Phi_{\alpha\beta}$ is the
fiber over $(\alpha,\beta)$. $\Phi$ is a principal
$\Gamma$-bundle, where $\Gamma=\{g\in\Aut(\pp^n): gY=Y\}$.If
$Y=X\cap\pp^{n-1}=X\cap H_{\alpha}$, consider the restriction
of $\Phi$ to $\{\alpha\}\times U$. Each loop $\{\alpha_t\}
(t\in [0;1])$ in $U$ can be lifted to this restriction as the
path $\{(\alpha_t;\varphi_t)\}$, where
$\varphi_t\in\Phi_{\alpha,\alpha_t}$. It is clear that
$\varphi_1$ induces in $H^\bullet(Y)$ the monodromy
transformation corresponding to the loop $\{\alpha_t\}$. The
proposition is proved.

\proclaim Corollary 1{.}2.
If a smooth variety $X\inPn$ satisfies the FF condition, then
the monodromy group of its hyperplane section is finite; if
$\dim X$ is even, this group is trivial.

{\bf Proof.} Since the connected component
$\Gamma_0\subseteq\Gamma$ acts trivially in cohomology and
$\Gamma/\Gamma_0$ is finite, the first assertion holds; the
second assertion follows immediately from the first one and the
Picard-Lefschetz theory.

The following proposition is quite similar to the main
construction of~[8], so we omit some details of the proof.

\proclaim Proposition 1{.}3.
Let $X\inPn, X\ne\pp^n$ be a  smooth  projective variety, and
$L\subseteq(\pp^n)^{\ast},L\cong\pp^1$
a Lefschetz pencil of hyperplanes.
Assume that there exist a Zariski open subset $U\subseteq L$
and $\alpha\in U$ such that for any $\beta\in U$
there exists a projective  isomorphism
$\map{\psi_\beta}{H_\beta}{H_\alpha}$
satisfying the following conditions:
\medskip
\item{i)} $\psi_\beta(H_\beta\cap X)=H_\alpha\cap X$;
\item{ii)}   $\psi_\beta$  is  identity  on
$H_\alpha\cap H_\beta$.
\medskip
Then $X$ is a quadric.

{\bf Proof.}
Choose the homogeneous coordinates in $\pp^n$ so that the
equations of $H_\alpha$ (resp.\ the axis of $L$) are $x_n=0$
(resp.\ $x_{n-1}=x_n=0$). For $u\in\cc$ denote by $H_u$ the
hyperplane defined by the equation $x_n=ux_{n-1}$, and denote
by $H_\infty$ the hyperplane $x_{n-1}=0$.

Now for $\beta\in U$ denote by $\Psi_\beta$ the set of
projective automorphisms $\map{\psi_\beta}{H_\beta}{H_\alpha}$,
such that $\psi_\beta(H_\beta\cap X)=H_\alpha\cap X$ and
$\Psi|_{H_\alpha\cap H_\beta\cap X}={\rm id}$. Consider a fiber
space $\Psi$ over $U$ such that $\Psi_\beta$ is the fiber over
$\beta$. Let $\Gamma\subseteq\Psi$ be a quasi-section of $\Psi$
over an open subset $U'\subseteq U$; the projection
$\map{\pi}{\Psi}{U}$ induces a regular function $u$ on $\Gamma$
such that any $p\in\Gamma$ may be regarded as a linear
isomorphism $f_p:H_u\to H_0=H_\alpha$; writing $f_p^{-1}$ in
the matrix form, we obtain regular functions
$a_0,\ldots,a_{n-1}$ such that $f_p^{-1}$ sends
$(x_0:\cdots:x_{n-1})\in H_{\alpha}$ to the point
$$
\eqalign{
(x_0+&a_0(p)x_{n-1}:\cdots:x_{n-2}+a_{n-2}(p)x_{n-1}:\cr
&a_{n-1}(p)x_{n-1}:u(p)a_{n-1}(p)x_{n-1})\in X\inPn.\cr}
$$
Set $a_{n}=u\cdot a_{n-1}$.
If $S$ is the smooth projective model of $\Gamma$, then
$a_j$'s may be regarded as rational functions on $S$; not
all of them are constant, because $u$ is not constant and
$a_{n}=u\cdot a_{n-1}$.

Now we can proceed as in~[8,Section 3]: not all
$a_j$'s, $0\le j\le n$, are constant, hence some of them must
have poles. Assume that the maximal order of these poles equals
$m$ and is attained at the point $\xi\in S$; by~[7, Lemma
3.1], for each $c\in\cc$ and $x=(x_0:\cdots:x_{n-1}:0)$
there exists a map $\map{h}{\Delta}{X\cap H_{\alpha}},
h:t\mapsto (\tilde{x}_0(t):\cdots:\tilde{x}_{n-1}(t):0)$, where
$\Delta$ is the unit disk in the complex plane, such that
$\tilde{x}_i(0)=x_i$ for $0\le i\le n-2$, $\tilde{x}_{n-1}\sim
ct^m$ as $t$ tends to 0. The point $\lim_{t\to
0}f_{h(t)^{-1}}(h(t))$ is in $X$, and its homogeneous
coordinates are
$$
  (x_0+cb_0x_{n-1}:\cdots:x_{n-2}+cb_{n-2}x_{n-1}:
     cb_{n-1}x_{n-1}:cb_nx_{n-1});
$$
as it is explained in~[8], $b_j$'s do not depend on $c$.
Hence $X\cap H_{b_n/b_{n-1}}$ is a cone;
since $L$ is a Lefschetz pencil,
this cone must be a quadratic cone and $X$ must be a quadric.

\proclaim Corollary 1{.}4.
If $X\inPn$ is not a linearly normal rational scroll nor the
Veronese surface $v_2(\pp^2)\subset\pp^5$, then Proposition~1.3
holds with hypothesis~(ii) replaced by {\it ``$\psi_\beta$ is
identity on $H_\alpha\cap H_\beta\cap X$''.}

{\bf Proof.}
If $\deg X=d$, then the hypothesis implies that its generic
linear section of codimension 2 contains at least $d+1$ points
in general position, hence ``identity on the linear section of
codimension 2'' implies ``identity on the projective space of
the section'', and the Proposition applies.

\bigskip
\centerline{\bf 2. Case of curves}\medskip
In this section we prove Proposition~0.1. Assume that $X\inPn$
is a curve for which FF holds, and that $X$ is not contained in
a hyperplane. We are to prove that $\deg X\le n+1$.

Let us apply Proposition~1.1. In our case the generic
hyperplane section is a set of $\deg X$ points in $\pp^{n-1}$,
and the monodromy group consists of permutations of these
points. According to [2], this group is the whole symmetric
group; on the other hand, if $Y=X\cap \pp^{n-1}$ is the generic
hyperplane section, then no $n$ points of $Y$ belong to a
hyperplane (we will call it the {\it generic position
property\/}). Proposition~0.1 follows immediately from the
above observations and the following lemma.
\proclaim Lemma 2{.}1.
If there are $s>m+2$ points in $\pp^m$ such that no $m+1$ of
them belong to a hyperplane, then there is no automorphism of
$\pp^m$ that interchanges two of these points and leaves the
rest $s-2$ points fixed.

{\bf Proof.} If $s\ge m+4$, there is nothing to prove since any
projective automorphism of $\pp^m$ fixing $m+2$ points in
general position must be identity. Hence we may assume that
$s=m+3$. Due to the generic position condition we may choose
the homogeneous coordinates so that $p_1=(1:0:\ldots:0),
p_2=(0:1:\ldots:0),\ldots, p_{m+1}=(0:\ldots:0:1),
p_{m+2}=(1:\ldots:1)$. If $p_{m+3}=(x_0:\ldots:x_m)$, then it
follows from the generic position condition that $x_i\ne 0$ for
all $i$, $x_i\ne x_j$ for $i\ne j$. Hence the automorphism
$\map{\varphi}{\pp^m}{\pp^m}$ that interchanges $p_{m+2}$ and
$p_{m+3}$ should be defined by a diagonal matrix
$\diag(x_0,\ldots,x_m)$; since $\varphi(x_{m+3})=x_{m+2}$, we
see that each of the $x_j$ can be chosen to equal 1 or $-1$;
this contradicts the generic position condition. The lemma and
Proposition~0.1 are proved.

\bigskip
\centerline{\bf  3. Case of surfaces, part 1}\medskip
We keep the notation of Section~1. Assume that $X\inPn$ is a
smooth surface for which the FF condition holds.
\proclaim Proposition 3{.}1.
If the generic hyperplane section of $X$ is not a rational
curve, then $\Phi\rightarrow U\times U$ is a finite covering.

{\bf Proof.} The fiber of $\Phi$ over
$(\alpha,\beta)\in U\times U$ is  isomorphic  to
$\{g\in\Aut H_{\alpha}:g(X\cap H_{\alpha})=X\cap H_{\beta}\}$.
Since the group of  automorphisms  of  a
smooth curve of genus $>1$ or a polarized elliptic curve is  finite,
we are done.

\proclaim Proposition 3{.}2.
If the genus of the generic hyperplane section of $X$ is
greater than 1, then the covering $\map{p}{\Phi}{U\times U}$ is
split.

{\bf Proof.} Assume the contrary; then the covering
$p^{-1}(U\times\{\alpha\})\rightarrow U\times\{\alpha\}$ is not
split for the generic $\alpha\in U$. Hence, there exists a
connected component $\Psi\subseteq p^{-1}(U\times\{\alpha \})$
such that $\map{p}{\Psi}{U\times\{\alpha\} }$ is a nontrivial
covering. Thus, there exists a loop in $U$ originating at
$(\alpha,\alpha)$, such that its lifting to $\Psi$ defines a
nontrivial automorphism of $X\cap H_{\alpha}$. Since any
nontrivial automorphism of a Riemann surface $C$ of genus $>1$
acts nontrivally on $H^1(C,\rr)$, we infer that this loop
defines a nontrivial element of the monodromy group of the
hyperplane section $X\cap H_{\alpha}$. This contradicts
Corollary~1.2.

\proclaim Proposition 3{.}3.
If $X\inPn$ is a smooth surface for which the FF condition
holds, then the genus of the generic hyperplane section of $X$
is at most 1.

{\bf Proof.} Assume the contrary, and let $\alpha\in U$ be a
generic point. Since $\Phi\rightarrow U\times U$ is split by
Proposition~3.2, there exists a section
$\map{s}{U\times\{\alpha\} }{\Phi}$, such that $s(
(\alpha,\alpha) )={\rm id}_{H_{\alpha} }$ Set
$\psi_{\beta}=\map{s(\beta)}{H_{\beta} }{H_{\alpha} }$. Define,
for the generic $x\in X$, the mapping $\map{f}{U}{X\cap
H_{\alpha} }$ by the formula $\beta\mapsto\psi_{\beta}(x)\in
H_{\alpha}$. By Proposition~1.4 this mapping is not constant,
hence $\overline{f(U)}=X\cap H_{\alpha}$. But this equality is
impossible since $X\cap H_{\alpha}$ is not a rational curve.
This contradiction completes the proof.

\bigskip
\centerline{\bf  4. Case of surfaces, part 2}\medskip

     In this section we assume  that $X\inPn$ is  a  linearly
normal smooth surface, that the condition FF holds for $X$
and  that  the
generic hyperplane section of $X$ is an elliptic curve. We make  use
of the following important result of Zak~[11]:

\proclaim  Theorem 4{.}1 (F{.}L{.}Zak).
If $X\inPn$ is a smooth surface such that the monodromy group
of hyperplane section of $X$ is trivial, then $X$ is either a
scroll, or the Veronese surface $v_2(\pp^2)$, or its isomorphic
projection to $\pp^4$.

It follows immediately from this theorem and Corollary~1.2 that
the assumptions of this section imply that $X$ is $\pp_C(\Ee)$
embedded by the complete linear system $|\Oo_{X|C}(1)|$, where
$\Ee$ is a rank 2 locally free sheaf on the elliptic curve $C$.
$H^0(\Oo_{X|C}(1))$ will be canonically identified with
$H^0(\Ee)$. Let us denote $\Ll=\det\Ee$.

If $s\in H^0(\Oo_{X|C}(1))=H^0(\Ee)$, consider the homomorphism
$\map{f_s}{\Ee}{\det\Ee}$ defined by the formula $\xi\mapsto
s\wedge \xi$.

\proclaim Proposition 4{.}1.
If the section $s$ defines  a  smooth  hyperplane
section of $X$, then the sequence of sheaves
$$
  \triple{\Oo_C}{s}{\Ee}{f_s}{\Ll}
\eqno (1)
$$
is exact.

The proof is straightforward.

\proclaim Proposition 4{.}2.
     Consider an exact sequence of sheaves
$$
  \triple{\Oo_C}{s}{\Ee}{f}{\Ll}
\eqno (2)
$$
where $C$ is an elliptic curve, $\Ee$ is a locally free sheaf of rank 2,
$\Ll$ is an invertible sheaf
(hence, $\Ll=\det\Ee$), as an extension of $\Ll$ by
$\Oo_C$. The class of this extension
in $\ext(\Ll,\Oo_C)$ is determined, up to
proportionality, by the linear subspace
$$
  \Im(H^0(C,\Ee)\rightarrow H^0(C,\Ll))\subseteq H^0(C,\Ll).
$$

{\bf Proof.} Since the sheaves are locally free and the
underlying variety is a smooth curve, $\ext(\Ll,\Oo_C)$ is
canonically isomorphic to $H^1(C,\Ll^{-1})$ and, by Serre's
duality, canonically dual to $H^0(C,\Ll)$; the fundamental
class of the extension~(2) in $\ext(\Ll,\Oo_C)=(H^0(C,\Ll)
)^\ast$ is $$\map{\delta}{H^0(\Ll)}{H^1(\Oo_C)}\cong\cc,$$
where $\delta$ is the connecting homomorphism associated with
the exact sequence~(2). Hence this class is determined, up to
proportionality, by $\Ker\delta=\Im f_\ast$.

Let us return to our surface $X$.
\proclaim Proposition 4{.}3.
For generic hyperplanes $H_1,H_2\subseteq\pp^n$ there exists a
projective automorphism $\map{F}{\pp^n}{\pp^n}$, such that
$F(X)=X$, $F(H_1)=H_2$, and $F$ maps each line of the ruling of
$X$ into itself.

{\bf Proof.} The smooth hyperplane section of $X$ defined by  a
section $s\in H^0(\Ee)$ is projectively isomorphic to
the curve $C$ embedded by
the linear  system $|V_s|$,  where
$$
  V_s=\Im(H^0(\Ee)\stackrel{(f_s)_\ast}{\rightarrow}H^0(\Ll) ).
$$
It follows immediately, from the exact cohomology sequence associated
with (1) and the ampleness of $\Ee$, that $V_s$ has
codimension $1$ in $H^0(\Ll)$.

Now if $s$ and $t$ are two generic sections of $\Ee$, then the
hyperplane sections defined by $s$ and $t$ are projectively
isomorphic if and only if there exists an isomorphism
$\map{\varphi}{C}{C}$ such that $\varphi^\ast\Ll=\Ll$ and
$V_s=\varphi^\ast V_t$. Since the group of automorphisms of a
polarized elliptic curve is finite, the FF condition implies
that the hyperplanes $V_s\subseteq H^0(\Ll)$ are the same for
almost all $s\in H^0(\Ee)$. By Proposition~4.2 this implies
that the extensions (1) are ``congruent up to multiplication by
a constant''for various $s$. Hence, for generic $s,t\in
H^0(\Ee)$ there exists an automorphism $\map{g}{\Ee}{\Ee}$ and
a constant $\lambda\in\cc^\ast$ such that the diagram
$$
  \matrix{
     0 & \longrightarrow & \Oo_C & \stackrel{s}{\longrightarrow} & \Ee
       & \stackrel{f_s}{\longrightarrow} & \Ll & \longrightarrow & 0\cr
     & & \Big\Vert & & \vertmap{g} & & \Big\Vert & & \cr
     0 & \longrightarrow & \Oo_C & \stackrel{\lambda t}{\longrightarrow} & \Ee
       & \stackrel{f_t}{\longrightarrow} & \Ll & \longrightarrow & 0\cr
  }
$$
is commutative. The automorphism $g$ induces a projective
automorphism $\map{F}{\pp^n}{\pp^n}$ that maps $X$ into itself
and preserves the fibers of $X$ over $C$. Translating all this
into the geometric language, we obtain our proposition.

\proclaim Proposition 4{.}4.
If $X\inPn$ is a linearly normal scroll with elliptic base,
then the FF condition does not hold for $X$.

{\bf Proof.} Assume the contrary. Then Proposition~4.3 applies.
Since generic hyperplane section intersects each line of the
ruling only once, the automorphism $F$ of the above proposition
fixes all the points of $H_1\cap H_2\cap X$. This contradicts
Proposition~0.2. The proposition is proved.

\bigskip
\centerline{\bf 5. Proof of propositions 0.2 and
0.3.}\medskip
To complete the proof of Proposition~0.2, we use the following
fact:

{\narrower\narrower\noindent
if the generic hyperplane section of a smooth  surface
face $X\inPn$ has genus 0, then $X$ is either a  rational  scroll,  or
$\pp^2$, or the Veronese surface
$v_2(\pp^2)$, or its  isomorphic  projection.

}
When put together with Propositions~3.3 and~3.4, this yields
the required result.

To prove Proposition~0.3, observe that if a scroll $X$ is
isomorphic to $\pp_C(\Ee)$, where $C$ is a curve and $\Ee$ is a
locally free sheaf of rank 2, such that $H$ is a hyperplane
section of $X$, then $(H,\Oo(1))\cong (C,\det\Ee)$, where $H$
is a hyperplane section of $X$. Hence, if this hyperplane
section were linearly normal, then almost all hyperplane
sections would be projectively isomorphic to the curve $C$
embedded by the complete linear system $|\det\Ee|$, contrary to
Proposition~0.2. The proposition is proved.

\bigskip
\centerline
{\bf 6. Appendix. Finite monodromy groups and adjunction}\medskip

The FF property and finiteness of the monodromy group have to do
with the adjunction properties of the variety.
\proclaim Proposition 6{.}1.
Consider the following properties of a smooth projective variety
$X\inPn, \dim X=d>1:$

\item{i)}   Almost all hyperplane sections of $X$ are projectively
equivalent.
\item{ii)}  The monodromy group of hyperplane sections of $X$ is finite.
\item{iii)}
If $p\ne q, p+q=d-1$,
then $h^{p,q}(X)=h^{p,q}(Y)$, where
$Y$ is a smooth hyperplane section of $X$.
\item{iv)} $|K_X+Y|=\emptyset$.

Then the following implications hold:
$$
  \hbox{i))}\Rightarrow\hbox{ii)}
  \Leftrightarrow\hbox{iii)}
  \Rightarrow\hbox{iv)}.
$$
If, in addition, $\dim X\le 3$, then
$\hbox{iii)}\Leftrightarrow\hbox{iv)}$.

{\bf Proof.} The implication $\hbox{i)}\Rightarrow\hbox{ii)}$
is just Corollary~1.2; the equivalence of ii) and iii) is
proved in~[3, Expos\'e XVIII]. To prove that
$\hbox{iii)}\Rightarrow\hbox{iv)}$, observe that $\hbox{iii)}$
implies that
$$
  h^{d-1}(X,\Oo_X)=h^{d-1}(Y,\Oo_Y).
\eqno (3)
$$
Now the exact sequence
$$
  \triple{\Oo_X(-1)}{}{\Oo_X}{}{\Oo_Y}
$$
together with Kodaira vanishing theorem, yields the exact sequence
$$
\eqalign{
  0\to &H^{d-1}(X,\Oo_X)\to H^{d-1}(Y,\Oo_Y)  \to \cr
       &  H^d(X,\Oo_X(-1))\to H^d(X,\Oo_X) \to 0.\cr
}
\eqno(4)
$$
Hence,(3) is equivalent to injectivity of the homomorphism
$H^d(X,\Oo_X(-1))\to H^d(X,\Oo_X)$ from~(4); by Serre
duality this is equivalent to the equality
$$
\dim|K_X|=\dim|K_X+Y|.
$$
The latter equality holds if and only if $|K_X+Y|=\emptyset$.
Indeed, the ``if'' part is obvious since $\dim|K_X|\le\dim|K_X+Y|$,
and to prove the ``only if'' part observe that $|K_X|\ne\emptyset$
implies the inequality $\dim|K_X+Y|>\dim|K_X|$, since the linear system
$|Y|$ is movable.

To prove the last assertion observe that, if $2\le\dim X\le 3$,
property~iii) is {\it equivalent\/} to the equality~(3). The
proposition is proved.

In the paper~[10], A.J.Sommese gave a complete description of
threefolds having property~iv). Proposition~6.1 shows that [10]
yields description of smooth threefolds with finite monodromy
group of hyperplane section, as well. All the threefolds with
FF property are among those from~[10]; no doubt only few of the
latter actually have the FF property.

\bigskip
\centerline{\bf  References}
\medskip
{
\frenchspacing\raggedright
\item{[1]} E.Ballico. {\it On projective varieties with
projectively equivalent zero-dimensional linear
sections,} Canad. J. Math. 35(1992), 3-13.

\item{[2]} E.Ballico, A.Hefez.
{\it On  the  Galios  group  associated  to  a
     generically \'etale morphism,} Commun. algebra. 14(1986),
     899--909.

\item{[3]} P.Deligne,  N.Katz.
{\it  Groupes  de  monodromie  en  g\'eom\'etrie
       alg\'ebrique,} Lect. notes in math., 340, Springer,
Berlin,   1973.

\item{[4]} G.Fano.
{\it Sulle superficie dello spazio $S_3$ a  sezioni  piane
      collineari, }
Rend R. Accad. dei Lincei, ser. $6^a$.
        1(1925),   473--477.

\item{[5]} G.Fano.
{\it  Sulle superficie di uno spazio qualunque  a  sezioni
      iperpiane collineari,} Rend  R.  Accad.  dei  Lincei,
      ser. $6^a$.
      2(1925),
      115--129.

\item{[6]} G.Fubini.
{\it  Sulle variet\`a a sezioni piane collineari,} Rend  R.
     Accad. dei Lincei,
      ser. $6^a$.
      1(1925),
      469--473.

\item{[7]} S.L'vovsky.
{\it Extensions   of   projective   varieties   and
     deformations I,} Mich. math. J, 39(1992), 41--51.

\item{[8]} S.L'vovsky.
{\it  Extensions   of   projective   varieties   and
     deformations II,}
Mich. math. J, 39(1992), 65--70.

\item{[9]} R.Pardini.
{\it Some remarks on projective varieties with projectively
equivalent hyperplane sections,}
preprint, Universit\`a di Pisa, fall~1992.

\item{[10]} A.J.Sommese.
{\it On the nonemptiness of the adjoint linear system
of a hyperplane section of a threefold,}
J. reine angew.\ Math., 402(1989), 211--220.

\item{[11]} F.L.Zak.
{\it Surfaces with zero  Lefschetz  cycles.} [in  Russian]
    Matematicheskiye zametki, 13(1973) 869--880.

}
\vfil\break
\obeylines
Serge L'vovsky
Russia, 117292, Moscow,
Profsoyuznaya ul., 20/9, kv.162.
\medskip%
nskcsmoscow@glas.apc.org
\bye